\documentclass[12pt]{amsart}
\usepackage{amsfonts}
\newcommand{\floor}[1]{\left\lfloor{#1}\right\rfloor}
\newcommand{\N}{{\mathbb N}}
\newtheorem{theorem}{Theorem}
\newtheorem{lemma}{Lemma}
\begin{document}

\title{Continued Fractions and Unique Additive Partitions}

\author{David J. Grabiner}\thanks{Research at MSRI supported by an NSF
Postdoctoral Fellowship.}

\email{grabiner@msri.org}

\address{Department of Mathematics, University of Michigan, Ann Arbor,
MI 48109-1003}

\begin{abstract}  A partition of the positive integers into sets $A$ and
$B$ {\em avoids} a set $S\subset\N$ if no two distinct elements in the
same part have a sum in $S$.  If the partition is unique, $S$ is {\em
uniquely avoidable.}  For any irrational $\alpha>1$, Chow and Long
constructed a partition which avoids the numerators of all convergents
to $\alpha$, and conjectured that the set $S_\alpha$ which this
partition avoided was uniquely avoidable.  We prove that the set of
numerators of convergents is uniquely avoidable if and only if the
continued fraction for $\alpha$ has infinitely many partial quotients
equal to 1.  We also construct the set $S_\alpha$ and show that it is
always uniquely avoidable.
\end{abstract}  

\maketitle

\section{Introduction}

A partition of the positive integers into sets $A$ and $B$ {\em avoids}
a set $S\subset\N$ if no two distinct elements in the same part have a
sum in $S$.  We say that $S$ is {\em avoidable}; if the partition is
unique, $S$ is {\em uniquely avoidable.}

The Fibonacci numbers are uniquely avoidable~\cite{aeh}.  Generalized
Fibonacci sequences defined by $(s_1,s_2)=1$, $s_{n}=s_{n-1}+s_{n-2}$,
are also uniquely avoidable provided that $s_1<s_2$ or $2|s_1s_2$;
Alladi, Erdos and Hoggatt~\cite{aeh} proved this for $s_1=1$, and
Evans~\cite{evans} proved the general case.  This suggests a connection
with continued fractions.  Chow and Long~\cite{chowlong} studied this
connection, and proved that the set of the numerators of
continued-fraction convergents to any irrational $\alpha$ with
$1<\alpha<2$ (easily generalized to any irrational $\alpha$) is
avoidable, although not necessarily uniquely avoidable.  Their partition
uses the sets
\begin{eqnarray*}
A_\alpha&=&\{n\in\N: \hbox{the integer multiple of $\alpha$ nearest $n$
is greater than $n$}\},\\
B_\alpha&=&\{n\in\N: \hbox{the integer multiple of $\alpha$ nearest $n$
is less than $n$}\}.
\end{eqnarray*}
Let $S_\alpha$ be the set avoided by the partition
$\{A_\alpha,B_\alpha\}$.  The main results of Chow and Long are that
$S_\alpha$ contains the numerators of all convergents to $\alpha$, and
every other element of $S_\alpha$ is either the numerator of an
intermediate fraction or twice the numerator of a convergent.

We give a characterization of a large class of sets which are uniquely
avoidable if they are avoidable at all.  We use this to show that the
set of numerators of convergents to $\alpha$ is uniquely avoidable for
any $\alpha$ if and only if infinitely many partial quotients of
$\alpha$ are 1.  We can also use the best approximation property to
determine $S_\alpha$ precisely, and show that $S_{\alpha}$ is uniquely
avoidable for any irrational $\alpha>1$.

\section{Results on continued fractions}

We will use some elementary results on continued fractions, as given
in~\cite{khinchin} for example.  We use the standard notation for
continued fractions, in which $\alpha=[a_0,a_1,\ldots]$, $p_{-2}=0$,
$q_{-2}=1$, $p_{-1}=1$, $q_{-1}=0$, and for $n\ge 0$, we have
$p_{n+1}=a_{n+1}p_n+p_{n-1}$.

For $1<\alpha<2$, let $\alpha'=\alpha/(\alpha-1)$.  Then
$\alpha=[1,a_1,a_2,\ldots]$, and $\alpha'=[a_1+1,a_2,\ldots]$.  Thus
$\alpha'$ and $\alpha$ have the same numerators of convergents;
$p_{0,\alpha'}=p_{1,\alpha}=a_1+1$, and the lost numerator
$p_{0,\alpha}=1$ appears as $p_{-1,\alpha'}$.  The numerators of
intermediate fractions between $p_{0,\alpha}=p_{-1,\alpha'}$ and
$p_{1,\alpha}=p_{0,\alpha'}$ are all the integers in this interval.
Also, since $1/\alpha + 1/\alpha'=1$, we have $A_{\alpha}=B_{\alpha'}$,
and thus $\alpha$ and $\alpha'$ give the same partition.  Thus
generalizing to arbitrary $\alpha$ does not give any extra generality;
we may assume $\alpha<2$ to reduce the number of special cases to be
considered.

We will study continued fractions in terms of their approximation
properties.  For any $p$, let $E(p)=p-q\alpha$ be the error in
approximating $p$ by the closest multiple of $\alpha$.  Then $A_\alpha$
is the set of all $n$ with $E(n)<0$.  Also, if $|E(x)+E(y)|<\alpha/2$,
then $E(x+y)=E(x)+E(y)$; otherwise, $E(x+y)=E(x)+E(y)\pm\alpha$.

We can restate several elementary results on continued fractions in
terms of $E$.  The best approximation property of convergents is that
$|E(x)|\le|E(y)|$ for all positive integers $y<x$ if and only if $x$ is
the numerator of a convergent.  The property that alternating
convergents approach $\alpha$ from opposite sides is that $E(p_n)<0$ for
$n$ even and $E(p_n)>0$ for $n$ odd.  The property that
$|\alpha-p_n/q_n|<1/(q_n q_{n+1})$ is that $|E(p_n)|<1/q_{n+1}$.

We also need the following lemma, which appears in similar form
in~\cite{chowlong}.
\begin{lemma}\label{doubleorint}
If $p<p_{n+1}$ and $|E(p)|<|E(p_{n-1})|$, then p is either $kp_n$ or
the numerator $kp_n+p_{n-1}$ of an intermediate fraction for some $k\in
N$. 
\end{lemma}

\begin{proof}
For simplicity of notation, we will assume $n$ is odd, so that
$E(p_n)>0$.  If $E(p)>0$, then let 
\[
k=\min\left(\floor{\frac{E(p)}{E(p_n)}},\floor{\frac{p}{p_n}}\right).
\]
Then $p'=p-kp_n$ still satisfies $0<E(p')<|E(p_{n-1})|$, and satisfies
either $p'<p_n$ or $0\le E(p')<E(p_n)$.  But no positive $p'<p_n$ can
have $0<E(p')<|E(p_{n-1})|$, by the best approximation property of
$p_{n-1}$, and no positive $p'<p_{n+1}$ can have $0<E(p')<E(p_n)$.  Thus
we have $p'=0$ and $p=kp_n$.

If $E(p)<0$, then $E(p-p_{n-1})=E(p)-E(p_{n-1})$, which is between $0$
and $-E(p_{n-1})$, so we can apply the above result to $p-p_{n-1}$.
\end{proof}

We will often need to use the following generalization of this lemma to
cases with $p>p_{n-1}$

\begin{lemma}\label{approxlemma}
If $|E(p)|<|E(p_{n-1})|$ and $p\le kp_{n+1}$, then $p=ip_{n-1}+jp_n$,
where $i$ and $j$ are non-negative integers and $i \le k$.
\end{lemma}

\begin{proof}
We again assume $n$ is odd, so that $E(p_{n+1})<0<E(p_n)$.  The proof is
by induction on $p$.  The previous lemma proves the case $p<p_{n+1}$.
If $E(p)>0$, then we let $p=p'+p_n$.  Since $E(p')<E(p)<|E(p_{n-1})|$
and $E(p')>-E(p_n)>-|E(p_{n-1})|$, and apply the lemma inductively to
$p'$ with the same $k$.  If $E(p)<0$, then we let
$p=p'+p_{n+1}=p'+p_{n-1}+a_{n+1}p_n$.  Since 
$E(p')>E(p)>-|E(p_{n-1})|$ and $E(p')<E(p_{n+1})<|E(p_{n-1})|$, we can
apply the lemma inductively to $p'$ with $k$ reduced by 1.
\end{proof}

\section{Characterization of $S_{\alpha}$}

The characterization of $A_\alpha$ and $B_\alpha$ in terms of $E$ gives
a natural characterization of the avoided set $S_\alpha$.

\begin{theorem}\label{testthm}
If $E(x)>0$, then $x\in S_\alpha$ if and only if there is no
even $z<2x$ with $0<E(z)<E(x)$; likewise, if
$E(x)<0$, then $x\in S_\alpha$ if and only if there is no
even $z<2x$ with $0>E(z)>E(x)$.
\end{theorem}

The following lemma, which we need in the proof of the theorem, is often
easier to use in showing that a particular $x$ does occur as a sum.

\begin{lemma}\label{testlemma}
If $E(x)>0$, then $x\in S_\alpha$ if and only if there is no $y<x$ with
$y\ne x/2$ and either $0<E(y)<E(x)$ or $0<E(y)+\alpha/2<E(x)$; likewise,
if $E(x)>0$, then $x\in S_\alpha$ if and only if there is no $y<x$ with
$y\ne x/2$ and either $0>E(y)>E(x)$ or $0>E(y)-\alpha/2>E(x)$.
\end{lemma}

\begin{proof}[Proof of Lemma~\ref{testlemma} and Theorem~\ref{testthm}]
If $E(x)>0$ and $x=y_1+y_2$ with $y_1,y_2\in B_\alpha$, and $y_1\ne
y_2$, then $E(y_1)$ and $E(y_2)$ are both positive and thus $E(y_1+y_2)$
is either negative or equal to $E(y_1)+E(y_2)$.  If $E(y_1+y_2)$ is
positive, then $E(y_1)$ and $E(y_2)$ are both less than $E(x)$.  Let
$E(y_1)$ be the smaller of the two (they are not equal since $y_1\ne
y_2$); then $E(y_1)<E(x)/2<\alpha/2$, and $0<E(2y_1)=2E(y_1)<E(x)$, so
we can take $z=2y_1$ in Theorem~\ref{testthm}.  Similarly, if
$x=y_1+y_2$ with $y_1,y_2\in A_\alpha$ and $y_1\ne y_2$, then $E(y_1)$
and $E(y_2)$ are both negative, and thus $E(y_1+y_2)$ can only be
positive if it is $E(y_1)+E(y_2)+\alpha$; it thus follows that
$E(y_1)+\alpha/2$ and $E(y_2)+\alpha/2$, which are both positive, must
be less than $E(x)$.  Let $E(y_1)+\alpha/2<E(x)/2$ be the smaller of the
two; then $z=2y_1$ gives $0<E(z)=2E(y_1)+\alpha<E(x)$.
\end{proof}

The following two theorems, the main theorems of~\cite{chowlong}, follow
immediately from Lemma~\ref{testlemma} and Theorem~\ref{testthm}.

\begin{theorem}\label{convergentthm}
The numerator $p_n$ of every convergent is in $S_\alpha$.
\end{theorem}

\begin{proof}
We apply Theorem~\ref{testthm} to $x=p_n$.  The number $z$ would have
to satisfy $|E(z)|<|E(p_n)|$, but Lemma~\ref{doubleorint} shows that any
such $z$ must be at least $p_n+p_{n+1}$, which is too large.
\end{proof}

\begin{theorem}\label{ambigthm}
If $p\in S_\alpha$, then $p$ is either the numerator $p_n$ of a
convergent, twice the numerator of a convergent, or the numerator
$p_{n-1}+kp_n$ of an intermediate fraction.
\end{theorem}

\begin{proof}
If $p$ is not the numerator of a convergent, then $p_n<p<p_{n+1}$ for
some $n$ and $|E(p)|>|E(p_n)|$.  As in Lemma~\ref{doubleorint}, we assume
that $n$ is odd for simplicity of notation.  If $0<E(p_n)<E(p)$, then we
can take $y=p_n$ in Lemma~\ref{testlemma} unless $p=2p_n$.  If
$E(p)<E(p_{n-1})<0$, then we can take $y=p_{n-1}$ in
Lemma~\ref{testlemma} unless $p=2p_{n-1}$.  Otherwise, we have
$E(p_n)<E(p)<-E(p_{n-1})$, which is only possible if $p$ is the numerator
of an intermediate fraction by Lemma~\ref{doubleorint}.
\end{proof}

The cases which are not resolved by these two theorems can be checked
similarly.

\begin{theorem}\label{doublethm}
For $p_n$ a convergent, $2p_n$ is in $S_\alpha$ if and only if $p_n$ is odd,
and either (i) $p_{n+1}$ is odd and $a_{n+1}\ge 3$, or (ii) $p_{n+1}$ is
even and $a_{n+1}\ge 2$, or (iii) $p_n=1$.
\end{theorem}

\begin{proof}
We first assume $|E(p_n)|<\alpha/4$; we will deal with the special case 
$|E(p_n)|>\alpha/4$ at the end.  In the general case, we will again
assume that $n$ is odd for simplicity of notation.

We will apply Theorem~\ref{testthm} to $x=2p_n$.  If $p_n$ is even,
then we can take $z=p_n$, and thus $2p_n\notin S_\alpha$.

If $p_n$ is odd, we need an even $z<4p_n$ with $0<E(z)<2E(p_n)$.  If
$E(z)<E(p_n)$, we can apply Lemma~\ref{approxlemma} to $z$.  If
$E(z)>E(p_n)$ with $z>p_n$, then we can apply Lemma~\ref{approxlemma} to
$z-p_n$.  If $E(z)>E(p_n)$ with $z<p_n$, then we can take $z'=2p_n-z$,
which has $0<E(z')<E(p_n)$, and use $z'$ instead of $z$.  Thus, in any
of these cases, we must have an even $z=ip_n+jp_{n+1}$ for some $i,j\ge
0$ with $z<4p_n$ and $0<E(z)<2E(p_n)$.  If $j$ is even, then since $p_n$
is odd, $i$ must also be even, but this is impossible because $z=2p_n$
is forbidden and $z=4p_n$ or $z=2p_n+2p_{n+1}$ is too large.  If $j=3$,
then we must have $i=0$ to have $z<4p_n$, but if
$|E(p_{n+1})|<\alpha/6$, we get $E(3p_{n+1})=3E(p_{n+1})<0$ while we
need $E(z)>0$, and if $\alpha/6<|E(p_{n+1})|<\alpha/4$, we have
$E(3p_{n+1})=\alpha+3E(p_{n+1})>\alpha/4$ while we need
$E(z)<E(p_n)<\alpha/4$.  If $j=1$, then we must have $i>0$ since
$E(p_{n+1})>0$.  The condition $z<4p_n$ is equivalent to $i<a_{n+1}-3$,
since $z=p_{n+1}+ip_n=p_{n-1}+a_{n+1}p_n+ip_n$, and $p_{n-1}<p_n$ for
$p_n\ne 1$.  Either choice of $z$ has $E(z)$ in the correct range, since
$E(ip_n+p_{n+1})=iE(p_n)+E(p_{n+1})$, which is positive because
$E(p_n)>-E(p_{n+1})$ and less then $2E(p_n)<\alpha/2$ since
$E(p_{n+1})<0$.  Thus we need only check that $z$ is even; it is even if
$p_{n+1}$ is odd and $i=1$, or $p_{n+1}$ is even and $i=2$.

We now need to check the special cases with $|E(p_n)|>\alpha/4$.  Recall
that we may assume $\alpha<2$.  Since $|E(p_n)|<1/q_{n+1}$, there can
only be problems if $q_{n+1}\le 3$, which implies $q_n\le 2$; thus the
only possible $p_n/q_n$ are $1/1$, $2/1$, and $3/2$.  Case (iii) covers
$p_n=1$.  If $p_n/q_n=2/1$ is a convergent, then $3/2<\alpha<2$; the
theorem says that the sum of 4 should occur, and it does occur because 1
and 3 are both in $A_\alpha$.  If $p_n/q_n=3/2$ is a convergent, then
$4/3<\alpha<2$, and for any such $\alpha$, $|E(p_n)|<\alpha/4$.
\end{proof}

\begin{theorem}\label{intthm}
The numerator $x=p_n+kp_{n+1}$ of an intermediate fraction is in
$S_\alpha$ if and only if either (i) $p_{n+1}$ is even, or (ii) $k=1$
and $p_n$ is odd, or (iii) $k=a_{n+2}-1$ and $p_{n+2}$ is odd.
\end{theorem}

\begin{proof}
Again, we will prove the special case $|E(p_n)|>\alpha/3$ at the end,
and assume in the general case that $n$ is odd for simplicity of notation.

Since $x$ is the numerator of an intermediate fraction,
$0<-E(p_{n+1})<E(x)<E(p_n)$.  We can apply Theorem~\ref{testthm} to see
whether $x\in S_\alpha$; it is not if there is an even $z<2x$ with
$0<E(z)<E(p_n)+kE(p_{n+1})$.  Since $0<E(z)<E(p_n)$, we can apply
Lemma~\ref{approxlemma} to write $z=ip_n+jp_{n+1}$; since $x<p_{n+2}$,
we have $z<2p_{n+2}$ and thus $i\le 2$.  If $i=0$ then we must have
$j\le 2k+1$ to have $z<2x$, and since $p_n>-a_{n+2}E(p_{n+1})$ and
$k<a_{n+2}$, we have $j|E(p_{n+1})|<2E(p_n)$.  If $E(p_n)<\alpha/4$,
then $E(z)=E(jp_{n+1})=jE(p_{n+1})<0$, which is impossible; if
$\alpha/4<E(p_n)<\alpha/3$, then we may instead have
$E(z)=\alpha+jE(p_{n+1})>\alpha-2E(p_n)>E(p_n)$, which is also
impossible.  If $i=2$, we have
$E(z)>E(z-p_{n+2})=E(z)-E(p_{n+2})=E(p_n+(j-a_{n+2})p_{n+1})$, and
$E(p_n+(j-a_{n+2})p_{n+1})$ is less than $E(p_n+kp_{n+1})$ only if
$j-a_{n+2}>k$, which implies $z>x+p_{n+2}$; this is impossible because
$x$ is a numerator of an intermediate fraction between $p_{n+1}/q_{n+1}$
and $p_{n+2}/q_{n+2}$.

Thus the only possibility is $i=1$, which gives $z=p_n+jp_{n+1}$.  We
have $0<E(z)<E(x)$ for $k<j\le a_{n+2}$, and $z<2x$ if
$p_n+jp_{n+1}<2p_n+2kp_{n+1}$, which requires $j\le 2k$ since
$p_n<p_{n+1}$.  If $p_{n+1}$ is even, then $p_n$ is odd, and thus no
choice of $j$ gives even $z$; this is case (i).  If $p_{n+1}$ is odd and
$k=1$, the only choice allowed is $j=2$, which gives odd $z$ if $p_n$ is
odd and even $z$ if $p_n$ is even; this is case (ii).  If $p_{n+1}$ is
odd and $k=a_{n+2}-1$, the only choice allowed is $j=a_{n+2}=k+1$, which
$z=p_{n+2}$; this is case (iii).  If $p_{n+1}$ is odd and we have any
other $k$, we can take either $j=k+1$ or $j=k+2$; one of these will give
$z$ even.

We now deal with the special cases with $|E(p_n)|>\alpha/3$, which
requires $q_{n+1}\le 2$ and thus $q_n=1$.  If $p_n/q_n=1/1$, then $n=1$,
and $|E(p_1)|>\alpha/3$ only for $\alpha>3/2$.  We must have $a_2=1$,
and thus $p_2=2$, so case (i) applies to all odd $x<2a_3+1$.  We have
$2-1/(a_3+1)<\alpha<2$.  Therefore every odd number up to $2a_3+1$ is in
$A_\alpha$ and every even number up to $2a_3$ is in $B_\alpha$, and thus
odd numbers less than $2a_3+1$ are avoided as required.  If
$p_n/q_n=2/1$. then $n=2$, $a_2=1$, and $\alpha>3/2$, so
$|E(2)|<1/2<\alpha/3$.
\end{proof}

Theorems \ref{convergentthm}, \ref{ambigthm}, \ref{doublethm}, and
\ref{intthm} give the complete characterization of the avoided set
$S_\alpha$.

\section{Uniquely avoidable sets}

The most natural characterization of unique avoidability is the
graph-theoretic characterization of~\cite{aeh}.  If $S\subset\N$, then
the {\em graph} $G(S)$ of $S$ is the graph with vertex set is $\N$ with
an edge between $x$ and $y$ if $x\ne y$ and $x+y\in S$.  A partition
which avoids $S$ is a 2-coloring of $G(S)$, and thus $S$ is avoidable if
and only if $G(S)$ is bipartite, and $S$ is uniquely avoidable if and
only if $G(S)$ is bipartite and connected.

Note that, if $a\in S$ and $x<a$, then $x$ and $a-x$ are in the same
connected component of $G(S)$; either $x=a-x$, or $x\ne a-x$ and their
sum is $a$.  Thus, if $a$ and $a+b\in S$, then $x$ and $x+b$ are in the
same connected component for every $x<a$.  This allows us to prove the
following result.

\begin{theorem}\label{uniquethm}
If $S$ contains $d$, $a+d$, $c$, and $b+c$ with $(a,b)=1$, $a\le c$, and
$b\le d$ (in particular, if $S$ contains $a$, $b$, and $a+b$), then all
numbers less than $a+b$ are in the same connected component of $G(S)$.
If $S$ contains infinitely many such subsets, then $G(S)$ is connected,
and therefore $S$ is uniquely avoidable if it is avoidable at all.
\end{theorem}

\begin{proof}
Assume $a<b$.  By the Chinese Remainder Theorem, any integer $m$ with
$0<m\le ab$ (and thus any $m$ with $0<m<a+b$) can be uniquely written
$m=xa-yb$ with $0<x\le b$, $0\le y<a$.  We prove connectivity for
$0<m<a+b$ by induction on $x+y$.  The base case is $m=a$.  For any other
$m$, if $m>a$, then $m$ and $m-a$ are in the same component; if $m<a$,
then $m$ and $m+b$ are in the same component.
\end{proof}

Our main results on unique avoidability follow.

\begin{theorem}
The set of numerators of convergents of $\alpha$ is uniquely avoidable
if and only if infinitely many partial quotients $a_n$ are 1.
\end{theorem}

\begin{proof}
If infinitely many $a_n$ are 1, then we can apply
Theorem~\ref{uniquethm} to $p_{n-2}$, $p_{n-1}$, and $p_n$ to show that
the graph is connected, and we have already shown that the set of
numerators is avoidable.  If only finitely many $a_n$ are 1, then there
is some $N$ such that for any $n>N$ we have $p_n>2p_{n-1}$.  Thus,
regardless of the partition of integers less than $p_{n-1}$, we can
extend the partition up to integers less than $p_n$ by placing integers
from $p_{n-1}$ to $\floor{p_n/2}$ arbitrarily, and placing $p_n-x$ in
the set which does not contain $x$ for $\floor{p_n/2}<x<p_n$.  We can
continue inductively for all $n$.
\end{proof}

\begin{theorem}
The set $S_\alpha$ avoided by the partition $A_\alpha,B_\alpha$ is
uniquely avoidable for any irrational $\alpha>1$.
\end{theorem}

\begin{proof}
If $p_{n+1}$ is even, then $p_n+p_{n+1}$ is either $p_{n+2}$ or the
numerator of an intermediate fraction, and if it is the numerator of an
intermediate fraction, it is in $S_\alpha$ by case (i) of
Theorem~\ref{intthm}.  If $p_{n+1}$ is odd, then $p_{n+1}-p_{n}$ is
either $p_{n-1}$ or the numerator of an intermediate fraction, and if it
is the numerator of an intermediate fraction, it is in $S_\alpha$ by
case (i) or (iii) of Theorem~\ref{intthm}.  In either case, we can apply
Theorem~\ref{uniquethm}.
\end{proof}

The unique avoidability of generalized Fibonacci sequences starting with
arbitrary $s_1,s_2$~\cite{evans} follows as a special case of these
results.  If $s_1<s_2$, we can let $s_1$ and $s_2$ be the numerators
$p_{n-1}$ and $p_n$ of two consecutive convergents, and let $a_m=1$ for
all $m>n$.  If $s_1>s_2$, we can let $p_n=s_2$, $p_{n+1}=s_1+s_2$; the
set of all the $s_i$ other than $s_1$ is thus uniquely avoidable, and
$s_1$ is the intermediate fraction $p_{n+1}-p_n$, which is an avoided
sum by case (i) of Theorem~\ref{intthm} if $s_2=p_n$ is even and by case
(iii) of Theorem~\ref{intthm} if $s_1+s_2=p_{n+1}$ is odd.  Thus the
sequence of $s_i$ is uniquely avoidable if $s_1<s_2$ or either $s_1$ or
$s_2$ is even, and not avoidable at all otherwise.


\begin{thebibliography}{9}
\bibitem{aeh} K. Alladi, P. Erd\H os, and V. E. Hoggatt, Jr., ``On additive
partitions of integers,'' {\em Discrete Math.} {\bf 23}(1978), 201-221.
\bibitem{chowlong} T. Y. Chow and C. D. Long, ``Additive partitions and
continued fractions,'' preprint.
\bibitem{evans} R. J. Evans, ``On additive partitions of sets of
positive integers,'' {\em Discrete Math.} {\bf 36}(1991), 239-245.
\bibitem{khinchin} A. I. Khinchin, {\em Continued Fractions}, 3rd ed.,
Chicago: University of Chicago Press, 1964.
\end{thebibliography}
\end{document}